\documentstyle[amstex,12 pt]{article}

\newtheorem{Th}{Theorem}[section]
\newtheorem{Prop}[Th]{Proposition}%[section]
\newtheorem{Lem}[Th]{Lemma}%[section]
\newenvironment{pf}{\noindent{\bf Proof.}}{\CQFD%\par\medskip 
}
\newcommand{\CQFD}
{%
\mbox{}%
\nolinebreak%
\hfill%
\rule{2mm}{2mm}%
\medbreak%
\par%
}

\newcommand{\R}{{\bf R}}

\newcommand{\C}{{\bf C}}
\newcommand{\D}{{\bf D}}

\def\Ho{\vbox{\offinterlineskip\hbox{\kern 3pt$\scriptstyle\circ$}
\kern
1pt\hbox{$H$}}}
\begin{document}

\title
{Locally Monge-Amp\`ere Parabolic Foliations}

\author{Morris KALKA and Giorgio PATRIZIO}
\date{ }

%\today
\footnotetext{Much of this work was done while Kalka was
visiting the University in Florence and G. Patrizio   Tulane
University. The Authors thank the instutions for their support. G. Patrizio
acknowledges the support of MIUR PRIN 2010-11 ``Variet\`a reali e complesse: geometria, topologia e analisi armonica'' and the collaboration with GNSAGA of INdAM.}

\maketitle
%\thispagestyle{empty}
%\newpage\addtocounter{page}{1}
\begin{abstract} 
\noindent
It is shown that codimension one parabolic foliations of complex manifolds are holomorphic.
This is proved using the fact that codimension one foliations of complex manifolds are necessarily locally Monge-Amp\`ere foliations and that  parabolic leaves cannot have hyperbolic behavior. 
The result holds true also for locally Monge-Amp\`ere foliations with parabolic leaves of arbitrary codimension.\\
\end{abstract}


\section {Introduction}

Let \(M\) be a complex manifold of complex dimension \(n>1\) and \(\cal F\) a complex foliation of \(M\) , i.e. a smooth foliation of \(M\) by complex submanifolds. We make the further assumption throughout that  \(\cal F\) is parabolic, i.e. the leaves of \(\cal F\) are parabolic complex manifolds  and of codimension \(p\) with $1\leq p< n$. By parabolic we mean that for any leaf \(L\) of \(\cal F\) there is a holomorphic cover
\[
F:\C^{n-p}\to L.
\]
By definition a a complex foliation \(\cal F\) of \(M\) is holomorphic if the distribution \({\mathcal T}\) tangent to \(\cal F\) is a holomorphic sub-bundle of \( T^{(0,1)}M\). 
 The geometry of parabolic foliations has been studied extensively in case the leaves are \(1\) dimensional and the foliation arises as the annihilator foliation of a plurisubharmonic solution of the complex homogeneous Monge-Amp\`ere equation
 \[
 (dd^cu)^n=0,
 \]
   under the nondegeneracy condition $ (dd^cu)^{n-1}\neq0$ (see for instance \cite{Bedford-Burns},\cite{Burns},\cite{Kalka-Patrizio2}, \cite{Kalka-Patrizio3},\cite{Patrizio},\cite{St}). We call such foliations, {\sl Monge-Amp\`ere foliations}.  In this situation one can use the parabolicity to prove uniformization type results for the manifold \(M\). Central to these results is the question of whether the foliation \(\cal F\) is holomorphic.   This question turns out to depend heavily on the global properties of the leaves, i.e. on parabolicity.

 In the codimension \(1\) situation, it is an easy result of Bedford-Kalka (\cite{Bedford-Kalka}) that any foliation by complex hypersurfaces  is {\it locally Monge-Amp\`ere} i.e. in the neighborhood of every point it is the annhilator foliation  of a plurisubharmonic solution of the complex homogeneous Monge-Amp\`ere equation
 \[
 (dd^cu)^{2}=0,
 \]
   under the nondegeneracy condition $ dd^cu\neq0$.
   
We find here that, using this fact, we can employ  techniques inspired by \cite{Bedford-Burns} and already used by Burns (\cite{Burns}) to prove  that any foliation by complex hypersurfaces,  whose leaves aresss uniformized by $\C^{n-1}$, is a holomorphic foliation.  

The technique for showing the holomorphicity of the foliation consists in proving the vanishing of a certain tensor. If \(\mathcal T\) denotes the complex tangent bundle to the leaves of \(\cal F\) and \(\mathcal N\) is the complex normal bundle to \(\cal F\), then  the Bedford-Burns twist tensor  
 ${\mathcal L}\colon{\mathcal T}\otimes \overline{{\mathcal N}}\to {\mathcal N}$ of the foliation is defined by
\begin{align}\label{twist}{\mathcal L}(V,\overline W)=[V,\overline W]\text{mod}({\mathcal T}\oplus T^{(0,1)}M),\end{align}
where the Lie bracket is computed using vector fields on $M$ which extend the  vector fields $V$ and $W$ defined along the leaves. It is known (see \cite {Bedford-Burns} and \cite  {Burns}  for example) that $\mathcal{L}$ vanishes on an open set $U$ if and only if the restriction of the foliation to $U$ is holomorphic. Indeed, the twist tensor defined by (\ref{twist}), commonly used as ``measure'' for (non-)holomorphicity in the the theory of Monge-Amp\`ere foliation, can be interpreted in terms of differential geometry of foliations. In particular it is known that ${\mathcal L}$  is equal to the antiholomorphic torsion of the complex Bott partial connection of the foliation. This is shown in \cite{Duchamp-Kalka} where these tools are used also to give necessary and sufficient conditions for a foliation in complex leaves to be locally Monge-Amp\`ere.
\\

For foliations of codimension $p>1$ there are examples (\cite{Duchamp-Kalka}, \cite{Calabi}) of non-holomorphic foliations whose leaves are all parabolic. These foliations, however, are shown  not to be locally Monge-Amp\`ere foliations. In fact we are able to show 
that, in any codimension, locally Monge-Amp\`ere foliations with parabolic leaves  are holomorphic. The idea of the proof is, again, an adaptation of arguments of Burns (\cite{Burns}).
In fact the hypersurface case (i.e. codimension $1$ foliation) is a special case of the general one.
We choose here to provide a separate proof for the codimension  $1$ case for which  the geometrical condition on the parabolicity of leaves is enough to prove the holomorphicity of the foliation and the curvature computations are much more transparent. The result for locally Monge-Amp\`ere foliations of arbitrary codimension is presented in the final section.

\section {Twist tensor and normal bundle  for locally Monge-Amp\`ere foliations of codimension 1}
 
We consider a codimension \(1\) locally Monge-Amp\`ere complex foliation \(\cal F\)  on a complex $n$-dimensional manifold $M$. We assume that \(\cal F\) is locally defined in a neighborhood $U$ of $q\in M$ by the annihilator 
 of the form $dd^{c}u$, where $u$ is a plurisubharmonic function on $U$ solving the complex homogeneous Monge-Amp\`ere equation
 \[
 (dd^cu)^2=0,
 \]
   under the nondegeneracy condition $dd^cu\neq 0.$
   
We will employ a local coordinate system on an open set \(U\) with respect to which the Bedford-Burns twist tensor takes a particularly convenient form and local computations for the normal bundle become easier. We refer to this coordinate system as {\sl leaf coordinates} and we define them as follows.
  If \(L\) is a leaf of the foliation  through $q$, we choose local coordinates $z_{1},\dots,z_{n}$ on an open set $U\ni p$, so that 
$L\cap U=\{z_{n}=0\}$, the functions $z_{1},\dots,z_{n-1}$ give holomorphic coordinates along the intersection $L\cap U$,   and, with the appropriate identification, $\frac{\partial \;\;}{\partial z_{n}}$ is a section on  the restriction to $L\cap U$ of the normal bundle to the foliation.
\begin{Prop}
Along $L\cap U$ the Levi matrix of $u$ has the form

\begin{align}\label{Levionleaf}\left({u_{j\bar k}}\right)_{|L\cap U} = \left(\begin{array}{cccc} 0 & \dots& \dots & 0
\\ \vdots & \dots & \dots & \vdots \\ \vdots  & \dots & 0 & 0 \\ 0  & \dots  & 0 & u_{n \bar n}    
\end{array}\right).
\end{align}
\end{Prop}
\begin{pf}
The Monge-Amp\`ere equation implies that $u$ is pluriharmonic along the leaf so that \(u_{j\overline k}=0\) for \(1\leq j,k\leq n-1\). Suppose  $W= \sum a^{i}\partial_{i}$ is a vector field tangent to a leaf. Then $W$ is in the annihilator of $dd^{c}u$ so that
$0=\sum \overline {a^{i}}u_{j\bar i}$ i.e. $\overline W(u_{j})= 0$ for all $j=1,\dots,n$ so that $u_{n}$ is holomorphic along the leaf i.e. $u_{n\bar i}=0$ for all $i=1,\dots, n-1$ in the leaf coordinates.
\end{pf}

We readily see that the Bedford-Burns twist tensor defined in (\ref{twist}) takes a particularly nice form in leaf coordinates. We consider a frame tangent to the leaves of the form ${\mathcal B}=\left\{  Z^{1},\dots, Z^{n-1}\right\}$ with  $Z^{j}=\frac{\partial\;}{\partial z_{j}}+b_{j}\frac{\partial\;}{\partial z_{n}}$.
The coefficients $b_{j}$ are determined by solving the system of equations $(0,\dots,1\dots,b_{j})(u_{l\bar k})=0$, i.e. by solving  
 the equation
 \[(0,\dots,1,0,\dots,0,b_{j})
 \left(\begin{array}{c}
 u_{1\,\overline{n}}\\
\vdots\\
u_{j\,\overline{n}}\\
\vdots\\
 u_{n\bar n}
 \end{array}\right)
 =0.
 \]
 for every $j$.
One immediately concludes that 
\begin{align}\label{foliationframe}  \hskip-.8cm
{\mathcal B}=\left\{  
Z^{1}=\frac{\partial\;}{\partial z_{1}}- \left(\frac{u_{1 \bar n}}{u_{n \bar n}}\right)\frac{\partial\;}{\partial z_{n}}, \dots
Z^{n-1}=\frac{\partial\;}{\partial z_{n-1}}- \left(\frac{u_{n-1 \bar n}}{u_{n \bar n}}\right)
\frac{\partial\;}{\partial z_{n}}
\right\}
\end{align}
defines a frame at every point of $U$ for the annihilator distribution $\mathcal A$.

\begin{Prop}
The components of the Bedford-Burns twist in leaf coordinates  are given  by:
\begin{align}\label{formulatwist}  
{\mathcal L}^{j}= 
\left[\left(u_{n \bar n}\right)^{-1}u_{\bar n  {j}}\right]_{\bar n}\frac{\partial\;}{\partial z_{n}}
\end{align}
for $j=1,\dots,n-1$, so that along the leaf $L\cap U$
\begin{align}\label{formulatwistleaf} {\mathcal L}^{j}=
\left(u_{n \bar n}\right)^{-1}u_{\bar n {j} \bar n}\frac{\partial\;}{\partial z_{n}}
\end{align}

 \end{Prop}
\begin{pf} 
Equality (\ref{formulatwist}) follows immediately
computing Lie brackets using the frame ${\mathcal B}$ defined by (\ref{foliationframe}), while (\ref{formulatwist}) is a direct consequence of (\ref{Levionleaf}).
\end{pf}

To make local computations regarding the normal bundle to the foliation, it is necessary to compute higher order derivatives of $u$. Under our assumptions we have the following

\begin{Lem} \label{forthderivatives} Let $L$ be a leaf of the Monge-Amp\`ere foliation and let $z_{1},\dots,z_{n}$ be leaf coordinates on $U$.
For $j,k\in  \{1,\dots,n-1\}$, if $A,B$ denote any pair of indices in the set $\{1,\dots,n, \bar1,\dots,\bar n\}$,
 one has :
\begin{align}\label{fourth} {u_{j\bar k AB}}_{|L\cap U}= {\left(u_{n\bar n}\right)^{-1}\left[u_{n\bar k A}u_{j\bar n B}+
u_{n\bar k B}u_{j\bar n A}\right]}_{|L\cap U}.
\end{align}
Moreover, unless $A,B\in \{n,\bar n\}$, one has ${u_{j\bar k AB}}_{|L\cap U}=0$.
 \end{Lem}
\begin{pf} With our choice of coordinates, along  the  intersection $L\cap U=\{z_{n}=0\}$ of the leaf  and $U$, the matrix $\left(u_{j\bar k}\right)$ is given by $(\ref{Levionleaf})$. Furthermore, the matrix $\left(u_{j\bar k}\right)$ has rank $1$ on $U$ because of our assumptions on the codimension of the foliation and hence $2 \times 2 $ minors of the matrix $\left(u_{j\bar k}\right)$ containing the term $u_{n\bar n}$ have determinant $0$.
Thus, for any $j,k\in  \{1,\dots,n-1\}$ and  any pair of indices $A,B$ in the set $\{1,\dots,n, \bar1,\dots,\bar n\}$, in order to compute $u_{j\bar k AB}$, we differentiate the equality:
\begin{align}\label{detzero}
0=\det
\left(\begin{array}{cc} 
u_{j \bar k}  &  u_{j \bar n} \\
u_{n \bar k}  & u_{n \bar n}    
\end{array}\right)= u_{j \bar k} u_{n \bar n} - u_{j \bar n} u_{n \bar k} 
\end{align}
and then use the particular form $(\ref{Levionleaf})$ that
 the matrix $\left(u_{j\bar k}\right)$ has on the leaf $L$. Thus from (\ref{detzero}) we have:
 \[
 \begin{array}{ll} 
 0 &= \left[ u_{j \bar k} u_{n \bar n} - u_{j \bar n} u_{n \bar k}   \right]_{AB}\\
 &\\
   &= \left[ u_{j \bar kA} u_{n \bar n} + u_{j \bar k} u_{n \bar nA}
   - u_{j \bar nA} u_{n \bar k} - u_{j \bar n} u_{n \bar kA}   \right]_{B}\\ 
   &\\
   & = u_{j \bar kAB} u_{n \bar n} + u_{j \bar kA} u_{n \bar nB}
   + u_{j \bar kB} u_{n \bar nA} + u_{j \bar k} u_{n \bar nAB}  \\
  & \\ 
 & \hskip2cm   - u_{j \bar nAB} u_{n \bar k} - u_{j \bar nA} u_{n \bar kB} 
   - u_{j \bar nB} u_{n \bar kA}  - u_{j \bar n} u_{n \bar kAB}.
\end{array}
\]
Along the leaf $L$ all terms $\left(u_{j\bar k}\right)$ with at least one index different from $n$ vanish and therefore so do  their derivatives with respect to vectors tangent to the leaf $L$. Hence the previous equality, on $L$ reduces to 

\begin{align} 
0= u_{n\bar n}u_{j\bar k AB}- u_{n\bar k A}u_{j\bar n B} - u_{n\bar k B}u_{j\bar n A}
\end{align}
which is (\ref{fourth}). Moreover, by the same argument, we get that
\[u_{j\bar k AB}= \left(u_{n\bar n}\right)^{-1}\left[u_{n\bar k A}u_{j\bar n B}+
u_{n\bar k B}u_{j\bar n A}\right]
=0\]
if at least one of the indices $A,B$ is not in $ \{n,\bar n\}$.
\end{pf}
\bigskip

For a   codimension \(1\) locally Monge-Amp\`ere complex foliation \(\cal F\)  on a complex $n$-dimensional manifold $M$ defined on an open set $U$  by the annihilator of the form $dd^{c}u$ where $u$ is a plurisubharmonic function on $U$ such that $
 (dd^cu)^2=0$ and  $ (dd^cu)\neq 0$, the form $dd^{c}u$ defines, in the obvious way, a metric on the normal bundle to the restriction  \(\cal F_{\rm U}\) to $U$ of the foliation. Using  Lemma \ref{forthderivatives}, one may compute the Ricci curvature $\phi$ of this metric on the normal bundle to the foliation defined by
$dd^{c}u$. Along a leaf $L$, with respect to leaf coordinates $z_{1},\dots,z_{n}$ on $U$, so that $z_{1},\dots,z_{n-1}$ are holomorphic coordinates along the intersection $L\cap U=\{z_{n}=0\}$ of the leaf  and $U$, one has
\begin{align}\label{formularicci}\phi= - \frac{i}{2}\sum_{j,k=1}^{n-1}[\log u_{n\bar n}]_{j,\bar k}dz_{j}\wedge dz_{\bar k}
\end{align}
\bigskip

 If the foliation is only locally Monge-Amp\`ere, it is not possible to define globally a metric on the normal bundle using the   information that it is locally defined by the annhilator of  the form $dd^{c}u$ for a plurisubharmonic function  $u$. In fact we only know that the manifold $M$ is covered by open sets  $U$ on which  \(\cal F\) is defined by the annihilator of the form $dd^{c}u$ where $u$ is a plurisubharmonic function on $U$ such that $ (dd^cu)^2=0$ and  $ (dd^cu)\neq 0$.
 On the other hand it is well known and easy to see (\cite{Bedford-Kalka}) that if the foliation is defined on the same open set by the annihilators of $dd^{c}u$ and $dd^{c}v$ for two different plurisubharmonic functions $u$ and $v$, then $dd^{c}u= \lambda dd^{c}v$ for some positive function $\lambda$ which is constant along the leaves of the foliation. Since \(\lambda\) is constant along the leaves of the foliation, the Ricci form (\ref{formularicci}) determined by \(u\) equals the Ricci form determined by \(v\). Therefore formula (\ref{formularicci}) 
 in fact  globally  defines a Ricci curvature on the normal bundle. 
 We can summarize 
 the conclusions of this discussion as follows:

\begin{Prop} Let  \(\cal F\) be a codimension \(1\) locally Monge-Amp\`ere complex foliation   on a complex $n$-dimensional manifold $M$. 
The Ricci curvature 
locally defined by  (\ref{formularicci}) defines globally  the Ricci curvature of the normal bundle to the foliation  \(\cal F\).
\end{Prop}

We have the following easy but important remark that relates the twist tensor to the Ricci curvature of the normal bundle to the foliation:

\begin{Prop} \label{ricciform} With respect to leaf coordinates $z_{1},\dots,z_{n}$ on $U$, the Ricci curvature  form defined in (\ref{formularicci}) is given by
$$\phi= - \frac{i}{2}\sum_{j,k=1}^{n-1}S^{j\bar k}
dz_{j}\wedge dz_{\bar k}=
- \frac{i}{2}\sum_{j,k=1}^{n-1}\left(u_{n\bar n}\right)^{-2}u_{\bar n j \bar n}u_{  n\bar k n}dz_{j}\wedge dz_{\bar k}$$
In particular for all $j=1,\dots,n-1$
\begin{align}\label{twistricci} 
0\leq S^{j\bar j}= \Vert {\cal L}^{j}\Vert ^{2}
\end{align}
where the norm is with respect to the metrics induced by 
$\left(u_{n\bar n}\right)$ and
$\left(u_{n\bar n}\right)^{-1}$ on the normal bundle to the foliation and its dual respectively.
Thus the foliation is holomorphic if and only if these terms vanish.
 \end{Prop}
\begin{pf} The proof reduces to a computation which uses the formulas for the fourth derivatives obtained in  Lemma \ref{forthderivatives}. In fact, with respect to leaf coordinates for  $L$
we have for $j, k\in \{1,\dots,n-1\}$:
\[
\begin{array}{ll}
S^{j\bar k} &= \left[\log u_{n\bar n}\right]_{j\bar k}= \left(u_{n\bar n}\right)^{-1}
        -\left(u_{n\bar n}\right)^{-2} u_{n\bar nj}u_{n\bar n\bar k} \\
        & \\
        &= \left(u_{n\bar n}\right)^{-2}\left[
        u_{n\bar k n}u_{j\bar n \bar n} + u_{n\bar k \bar n}u_{j\bar n n}
        -u_{n\bar nj}u_{n\bar n\bar k}
       \right] \\
       & \\
        &= \left(u_{n\bar n}\right)^{-2}u_{\bar n j \bar n}u_{  n\bar k n}.
\end{array}
\]
\end{pf}

\section {Codimension 1 foliations with parabolic leaves are holomorphic}

The main result on foliations whose leaves are complex parabolic hypersurfaces is the following:

\begin{Th} \label{codimensionone} A  complex foliation of rank $n-1$ with parabolic leaves on a $n$-dimensional complex manifold is holomorphic. 
 \end{Th}
 
\begin{pf} Suppose that ${\cal F}$ is a foliation of complex codimension $1$ on the complex $n$-dimensional manifold M.   As already remarked, it is well known (\cite{Bedford-Kalka}) that the foliation ${\cal F}$ is locally defined in a neighborhood U of any $q \in M$ by the annihilator  of the form $dd^{c}u$ where u is a plurisubharmonic function on U. Therefore we can fully take advantage of the computations performed in Section $2$.
Let $L$ be a leaf of the Monge-Amp\'ere foliation. We need to show that  the 
twist tensor ${\cal L}$ vanishes at all points $q\in L$. Since, by hypothesis $L$ is parabolic,  there is a holomorphic covering map $F\colon \C^{n-1} \to L$. For any $q\in L$, up to reparametrization, one may assume that $q=F(0)$. Let $I_{j}\colon \C\to \C^{n-1}$ be the $j$-th canonical injection: $I_{j}(z)=(0,\dots,z,\dots,0)$ and $f_{j}=F\circ I_{j} \colon \C\to  L$.
Then,  $L^{j}=f_{j}(\C)$ is a parabolic curve in $M$   for all $j=1,\dots,n-1$ with $q=f^{j}(0)\in L^{j}$. For all $j$ let $\iota_{j}\colon L^{j}\to L$ be the injection and denote with $\phi^{j} = \iota_{j}^{*}\phi$ the pull-back of the Ricci form $\phi$ of the metric on the normal bundle to $L$ defined by the form $dd^{c}u$. Then the form 
$$\psi^{j}=-\phi^{j}= \frac{i}{2} S^{j\bar j}dz_{j}\wedge dz_{\bar j}$$
is non-negative on $L^{j}$ and hence it defines a (pseudo$-$)metric on it. According to (\ref{twistricci}),
we need to show that the function $S=S^{j\bar j}$ vanishes identically along $L^{j}$.
We have the following:

\begin{Lem} \label{lemmaricci}
With the notation as above, if $\psi^{j}(q)>0$, then at $q$, for each $j$
\begin{align}\label{ricci}
Ric(\psi^{j})=-2\psi^{j} 
\end{align}
\end{Lem} 

\bigskip
\noindent
We postpone the proof of the Lemma until the end of the paragraph and we complete our argument.
The equality shown in the Lemma \ref{lemmaricci} implies that one can define a metric of negative constant curvature on the part of the curve $L^{j}$ where $\psi^{j}\neq 0$. By an Ahlfors Lemma argument, we now show that this cannot happen on a parabolic Riemann surface. For $R>0$, let 
$\D(R)=\{z\in \C \mid \vert z  \vert < R\}$. Then if we consider  $f_{j}\colon \D(R) \to L^{j}$, i.e. the restriction of the map $f_{j}$ to the disk $\D(R)$, then $\omega= f_{j}^{*}(\psi^{j})$ satisfies 
$Ric(\omega)=-2\omega$ because of (\ref{ricci}), that is to say that the metric defined by 
$\omega$ on the part of $\D(R)$ where $\omega\neq0$ is of constant curvature $-1$. Because of Ahlfors' Lemma, $\omega$ is dominated by the hyperbolic metric of $\D(R)$ defined by
$$\omega_{0}= \frac{i}{2} \frac{R^{2}}{ (R^{2}-\vert z\vert^{2})^{2}}$$  
i.e. $\omega\leq \omega_{0}$. If we write $\psi^{j}$ as
$\psi^{j}= \frac{i}{2} Sdz_{j}\wedge dz_{\bar j}$, this
 implies   that the function $S$ is such that
\begin{align}\label{estimate} S(f_{j}(0))\vert {f_{j}}'(0)\vert^{2}\leq \frac{1}{R^{2}}.\end{align}
As $f_{j}$ is a covering map, $\left({f_{j}}\right)'(0)\neq0$ so that, since (\ref{estimate}) holds for all $R>0$,
it follows that there is a contraddiction if $S(f_{j}(0))=S^{j\bar j}(p)\neq0$. 
In all these considerations $j=1,\dots,n-1$ was arbitrary and 
$q$ was any point in the leaf $L$.  By Lemma  3.2,  it follows that the twist tensor $\cal L$ vanishes at  any $q\in L$. Since $L$ is any leaf, the foliation is holomorphic.
\end{pf}
\bigskip
Finally, here is the proof we still need to provide: \\  \\
\noindent
{\bf Proof of Lemma \ref{lemmaricci}.}  
Since $F$ is a holomorphic covering map, then $dF(0)\neq0$. Hence
$F$ provides local coordinates $z_{1},\dots,z_{n-1}$ along the leaf $L$ in a neighborhood (on $L$) of  $q$. We extend these coordinates to holomorphic coordinates $z_{1},\dots,z_{n-1}, z_{n}$ for $M$ in a neighborhood of $q$ so that, locally, $L$ is given by $z_{n}=0$, i.e. leaf coordinates. With respect to these coordinates, the point $q$ is the origin. Furthermore,  the computations of Proposition \ref{ricciform}  show that 
$$ \psi^{j}= \frac{i}{2} S^{j\bar j}dz_{j}\wedge dz_{\bar j}=
\frac{i}{2} \left(u_{n\bar n}\right)^{-2}\left| u_{\bar n j \bar n}\right|^{2}dz_{j}\wedge dz_{\bar j}.$$
For simplicity, we denote 
$$S= S^{j\bar j}= \left(u_{n\bar n}\right)^{-2}\left| u_{\bar n j \bar n}\right|^{2}.$$
Proving our claim is equivalent to showing that
$$ \left[\log S\right]_{j\bar j}(0)=2S(0).$$
Here and thereafter, we remind that lower indices for functions denote derivatives in the usual manner. Using (\ref{fourth}) with $k=j$ and $A=B=n$, we have:
$$S_{j}= \left[ \left(u_{n\bar n}\right)^{-2}\left| u_{\bar n j \bar n}\right|^{2}   \right]_{j} \hskip8.9cm
$$
$$
=-2 \left(u_{n\bar n}\right)^{-3} u_{n\bar nj}\left| u_{\bar n j \bar n}\right|^{2} +
\left(u_{n\bar n}\right)^{-2} u_{\bar n j \bar n j} u_{\ n \bar  j n}+
\left(u_{n\bar n}\right)^{-2} u_{\bar n j \bar n } u_{n \bar  j n j}
$$
$$
=-2\left(u_{n\bar n}\right)^{-3} u_{n\bar nj}\left| u_{\bar n j \bar n}\right|^{2} +
\left(u_{n\bar n}\right)^{-2} u_{\bar n j \bar n j} u_{n \bar  j n}  \hskip3.7cm
$$
$$ \hskip2.7cm
+\left(u_{n\bar n}\right)^{-3}   u_{\bar n j \bar n } 
\left(u_{n \bar j  n }u_{ j \bar n n }  +     u_{n \bar j  n } u_{ j \bar n n } \right)
$$
$$
= \left(u_{n\bar n}\right)^{-2} u_{\bar n j \bar n j} u_{n \bar  jn}. \hskip8.2cm
$$
Differentiating again, we have
\begin{align}\label{Sjbarj}
\begin{array}{ll}
S_{j\bar j}&= \left[\left(u_{n\bar n}\right)^{-2} u_{\bar n j \bar n j} u_{n \bar  jn}   \right]_{\bar j} \\
&\\
&= -2\left(u_{n\bar n}\right)^{-3} u_{n\bar n \bar j} u_{\bar n j \bar n j} u_{n \bar  j n} 
+\left(u_{n\bar n}\right)^{-2} u_{\bar n j \bar n j \bar j} u_{n \bar  j n} 
\\
&\\
&\hskip2cm  +\left(u_{n\bar n}\right)^{-2} u_{\bar n j \bar n j} u_{n \bar  j n \bar j}.
\end{array}
\end{align}
To proceed, we like to eliminate the fifth derivative term using formula (\ref{fourth}) of Lemma \ref{forthderivatives} in the case $k=j$ and $A=B= \bar n$:

\[
\begin{array}{ll}
u_{\bar n j \bar n j \bar j}  &= \left[ u_{ j \bar j\bar n  \bar n}\right]_{j}
=\left[  2 \left(u_{n\bar n}\right)^{-1}u_{n\bar j \bar n}u_{j\bar n \bar n} \right]_{j} \\
&\\
& =  2\left[-\left(u_{n\bar n}\right)^{-2}u_{n\bar n j}u_{n\bar j \bar n}u_{j\bar n \bar n}
+  \left(u_{n\bar n}\right)^{-1}
 u_{n\bar j \bar n j}u_{j\bar n \bar n}
+  \left(u_{n\bar n}\right)^{-1} u_{n\bar j \bar n}u_{j\bar n \bar nj}\right]
\\
&\\
& =  -2\left(u_{n\bar n}\right)^{-2}u_{n\bar n j}u_{n\bar j \bar n}u_{j\bar n \bar n} 
+ 2 \left(u_{n\bar n}\right)^{-1} u_{n\bar j \bar n}u_{j\bar n \bar nj}\\
&\\ 
 & \hskip1cm +2 \left(u_{n\bar n}\right)^{-2}   u_{j\bar n \bar n} 
 \left[   u_{n\bar j n}u_{j\bar n \bar n}+  u_{n\bar j \bar n}u_{j\bar n n} \right]\\
 &\\ 
 &= 2 \left(u_{n\bar n}\right)^{-1} u_{n\bar j \bar n}u_{j\bar n \bar nj}
 + 2 \left(u_{n\bar n}\right)^{-2}   u_{j\bar n \bar n}  u_{n\bar j n}u_{j\bar n \bar n}.
\end{array}
\]
We now plug this expression for the fifth derivatives in (\ref{Sjbarj}):
\[
\begin{array}{ll}
S_{j\bar j}&= 
 -2\left(u_{n\bar n}\right)^{-3} u_{n\bar n \bar j} u_{\bar n j \bar n j} u_{n \bar  j n} 
+\left(u_{n\bar n}\right)^{-2} u_{\bar n j \bar n j \bar j} u_{n \bar  j n} 
+\left(u_{n\bar n}\right)^{-2} u_{\bar n j \bar n j} u_{n \bar  j n \bar j}\\
&\\
&= -2\left(u_{n\bar n}\right)^{-3} u_{n\bar n \bar j} u_{\bar n j \bar n j} u_{n \bar  j n} 
\\
&\\
 &\hskip2cm +2 \left(u_{n\bar n}\right)^{-3} u_{n\bar j \bar n}u_{j\bar n \bar nj}
+ 2 \left(u_{n\bar n}\right)^{-4}   u_{j\bar n \bar n}  u_{n\bar j n}u_{j\bar n \bar n}
u_{n \bar  j n} \\
&\\
& \hskip2cm +\left(u_{n\bar n}\right)^{-2} u_{\bar n j \bar n j} u_{n \bar  j n \bar j}\\
&\\
&= \left(u_{n\bar n}\right)^{-2} u_{\bar n j \bar n j} u_{n \bar j n \bar j} +
2 \left(u_{n\bar n}\right)^{-4}  \left(u_{j\bar n\bar n }\right)^{2} \left(u_{\bar j n  n}\right)^{2}\\
&\\
&= \left(u_{n\bar n}\right)^{-2} \left| u_{\bar n j \bar n j}\right|^{2}+
2 \left(u_{n\bar n}\right)^{-4}  \left|u_{j\bar n\bar n }\right| ^{4}.
\end{array}
\]

It is not restrictive to assume that $u_{n\bar n}(0)=1$. Thus:
$$S(0)= \left| u_{\bar n j \bar n}(0)\right|^{2},$$
$$\left|S_{j}(0)\right|^{2}=\left| u_{\bar n j \bar n}(0)\right|^{2}\left|u_{\bar n j \bar n j}(0)\right|^{2}$$
$$
S_{j\bar j}(0)=   \left| u_{\bar n j \bar n j}(0)\right|^{2}+
2  \left|u_{j\bar n\bar n }(0)\right| ^{4} 
$$
so that
$$ S(0)\left[\log S\right]_{j\bar j}(0)=
S_{j\bar j}(0)- \frac{\left|S_{j}(0)\right|^{2}}{S(0)}= 2  \left|u_{j\bar n\bar n }(0)\right| ^{4} = 2(S(0))^{2}
$$
which is equivalent to the equality we were seeking.
\CQFD

\section {Higher codimensional parabolic foliations}

In this section we extend the previous results to complex foliations of general codimension. Namely, we consider the following situation. Let \(M\) be a complex manifold of complex dimension \(n>1\) and \(\cal F\) a complex foliation of \(M\). We assume throughout that  \(\cal F\) is parabolic and of codimension \(p\), with $1\leq p< n$. By parabolic we mean that if \(L\) is a leaf of \(\cal F\) there is a holomorphic cover
\[
F:\C^{n-p}\to L.
\]
Furthemore we assume that \(\cal F\) is a {\sl locally Monge-Amp\`ere foliation} i.e.
for any given point $q\in M$ there exist an open set $U\ni q$ and  a plurisubharmonic function  $u\colon U\to \R$  satisfying
\begin{align}\label{pMA}
 (dd^cu)^{p+1}=0 \hskip1cm{\rm with} \hskip1cm (dd^cu)^{p}\neq0
\end{align}
   and
such that 
\({\cal F}_{|U}={\rm Ann} ({\rm dd^{c}}u)\).

 It is easy to see (see Lemma $3.1$ of \cite{Bedford-Kalka} and its proof) that if $v$ is another plurisubharmonic function defined on $U$ such that \({\cal F}_{|U}={\rm Ann} ({\rm dd^{c}}v)\), 
there exists a positive function $\lambda$ which is constant along the intersections of the leaves of \(\cal F\) with $U$ and such that 
\begin{align}\label{annihilator}
(dd^{c}u)^{p}=\lambda (dd^{c}v)^{p}.
\end{align}
Also in this setting it is crucial for computation to consider appropriate systems of coordinates which we call again {\sl leaf coordinates}. They are defined as follows.
If \(L\) is a leaf of the foliation  through $q$, one may choose local coordinates $z_{1},\dots,z_{n}$ on an open set $U\ni p$, so that 
$L\cap U=\{z_{n-p+1}=\dots,z_{n}=0\}$, 
the functions $z_{1},\dots,z_{n-p}$ give holomorphic coordinates along the intersection $L\cap U$,   and 
$\frac{\partial   \; \;\; \; \;\;}{\partial z_{n-p+1}},\dots,\frac{\partial  \;\;}{\partial z_{n}}$, with the appropriate identification, is a section on  the restriction to $L\cap U$ of the normal bundle to the foliation.

\begin{Prop} Suppose that on $U$ the foliation \({\cal F}\) is defined as the annihilator of a plurisubharmonic function $u$:  \({\cal F}_{|U}={\rm Ann} ({\rm dd^{c}}u)\). Then,
along $L\cap U$, the Levi matrix of $u$ is as follows:
\begin{align}\label{Levionleafhigher}\left({u_{j\bar k}}\right)_{|L\cap U} =
\left(\begin{array}{cccccc}
 0 & \dots& 0 & 0& \dots& 0\\
  \vdots &  \vdots& \vdots&\dots & \dots & \vdots \\ 
  \vdots &  \vdots& 0&0 & \dots &0 \\ 
  0  & \dots  & 0 & u_{{n-p+1}\, \overline{n-p+1}}& \dots& u_{{n-p+1}\, \bar n}\\
\vdots &  \vdots& \vdots&\vdots & \vdots & \vdots \\
  0  & \dots  & 0 & u_{n\, \overline{n-p+1}}& \dots& u_{n \bar n}    
\end{array}\right).
\end{align}
\end{Prop}
\begin{pf}
The Monge-Amp\`ere equation (\ref{pMA}) implies that $u$ is pluriharmonic along the leaf so that \(u_{j\overline k}=0\) for \(1\leq j,k\leq n-p\). Now suppose  $W= \sum a^{i}\partial_{i}$ is a vector field tangent to a leaf. Then $W$ is in the annihilator of $dd^{c}u$ so that
$0=\sum  \overline {a^{i}}u_{j\bar i}$ i.e. $\overline W(u_{j})= 0$ for all $j$ so that $u_{k}$, for 
$k=n-p+1,\dots,\leq n$, is holomorphic along the leaf i.e. $u_{k\bar i}=0$ for all $k=n-p+1,\dots,\leq n$ and $i=1,\dots, n-p$ in the leaf coordinates.
\end{pf}
       
Also in the higher codimension case, for locally Monge-Amp\`ere foliations the Bedford-Burns twist tensor  defined in (\ref{twist}) has a nice expression in leaf coordinates. Suppose that on $U$ the foliation \({\cal F}\) is defined as the annihilator of a plurisubharmonic function $u$:  \({\cal F}_{|U}={\rm Ann}( {\rm dd^{c}}u)\). If \(L\) is a leaf of the foliation  through $q$, choose leaf  coordinates $z_{1},\dots,z_{n}$ on  $U\ni q$ -- here, if necessary we shrink $U$, so that 
$L\cap U=\{z_{n-p+1}=\dots,z_{n}=0\}$, 
the functions $z_{1},\dots,z_{n-p}$ give holomorphic coordinates along the intersection $L\cap U$,   and 
$\frac{\partial \;\;\;\;\;\;\;}{\partial z_{n-p+1}},\dots,\frac{\partial \;}{\partial z_{n}}$, with the appropriate identification, is a section of  the restriction to $L\cap U$ of the normal bundle to the foliation. 
In these coordinates, a frame tangent to the foliation has the form:
\[
{\mathcal B}=\left(Z^{1}=\frac{\partial\;}{\partial z_{1}}+\sum_{l=n-p+1}^{n}b_{jl}\frac{\partial\;}{\partial z_{l}},\dots, Z^{n-p}=\frac{\partial\;\;\;\;}{\partial z_{n-p}}+\sum_{l=n-p+1}^{n}b_{n-p\,l}\frac{\partial\;}{\partial z_{l}}\right).
\] 
Denote $B= \left(b_{jl}\right)$
and
\begin{align}\label{minors}
 \hskip-1.5cm\Lambda=
\left( \begin{array}{ccc}u_{1\,\overline{n-p+1}}&\dots &u_{1\bar n}\\
\vdots&&\vdots \\
u_{n-p\,\overline{n-p+1}}&\dots& u_{n-p\bar n}
\end{array} \right), \hskip0.2cm
H=
\left( \begin{array}{ccc}
u_{n-p+1\,\overline{n-p+1}}&\dots &u_{1\,\bar n}\\
\vdots&&\vdots \\
u_{n\,\overline{n-p+1}}&\dots &u_{n\,\bar n}
\end{array}\right).
\end{align}
Notice that $H$ is an invertible matrix because of the non-degeneracy condition prescribed in (\ref{pMA}). Since each $Z^{j}$ must be in the annihilator of $dd^{c}u$, from the system of equations 
 \[(0,\dots,1\dots,0,b_{j\,n-p+1,\dots},b_{jn})
 \left(\begin{array}{ccc}
 u_{1\,\overline{n-p+1}}&\dots u_{1\bar n}\\
\vdots&\vdots&\vdots \\
u_{n-p\,\overline{n-p+1}}&\dots u_{n-p\bar n}\\
u_{n-p+1\,\overline{n-p+1}}&\dots u_{1\bar n}\\
\vdots&\vdots&\vdots \\
u_{n-p+1\,\overline{n-p+1}}&\dots u_{n-p\bar n}
 \end{array}\right)
 =\left(\begin{array}{c}0\\ \vdots \\0 \end{array}\right)
 \]
 one concludes that 
 $B= -\Lambda H^{-1}$
so that we may write:
\begin{align}\label{foliationframehigher}\begin{array}{ll}
{\mathcal B}&=\left(Z^{1}=\displaystyle{\frac{\partial\;}{\partial z_{1}}}+\displaystyle{\sum_{l=n-p+1}^{n}}b_{jl}\frac{\partial\;}{\partial z_{l}},\dots, Z^{n-p}=\frac{\partial\;\;\;\;\;}{\partial z_{n-p}}+\displaystyle{\sum_{l=n-p+1}^{n}}b_{n-p\,l}\frac{\partial\;}{\partial z_{l}}\right)\\
&\\
&=\left(\frac{\partial\;\;}{\partial z_{1}},\dots,\frac{\partial\;\;\;\;}{\partial z_{n-p}},
\frac{\partial\;\;\;\;\;\;\;\;\;}{\partial z_{n-p+1}},
\dots,\frac{\partial\;\;}{\partial z_{n}}\right)
\left(\begin{array}{c}
I_{n-p}\\-\Lambda H^{-1}
\end{array}\right).
\end{array}
\end{align}

Finally we recall that, in  leaf coordinates,  the Bedford-Burns twist is
given by 
\[{\mathcal L}=\left({\mathcal L}^{1}=\sum_{m=n-p+1}^{n}{\mathcal L}^{1}_{m},\dots,
{\mathcal L}^{n-p}= \sum^{n}_{m=n-p+1}{\mathcal L}^{p}_{m}\right),\]
where, for $j=1,\dots,n-p$ and  $m=n-p+1,\dots, n$
we have:
\begin{align} 
{\mathcal L}^{j}_{m}=\left({\mathcal L}(Z^{j},\frac{\partial\;}{\partial \bar z_{m}})\right)
=\left([Z^{j},\frac{\partial\;}{\partial \bar z_{m}}]~mod({\mathcal T}\oplus T^{(0,1)}M) \right).
\end{align}
Using the notations introduced in (\ref{minors}), 
and denoting  from now on 
\begin{align} \label{inverseH}
H^{-1}= {\cal H}= ({\cal H}^{r \bar l}).
\end{align}
the inverse matrix of $H$ in (\ref{minors}), we are ready to provide the expression in leaf coordinates of the Bedford-Burns twist:
\bigskip
\begin{Prop}
The components of the Bedford-Burns twist with respect to the leaf coordinates are given  by:
\begin{align}
{\mathcal L}^{j}_{m}= \sum_{l=n-p+1}^{n} \left[\sum_{r}(-\Lambda_{j r}){\cal H}^{r \bar l}\right]_{\bar m}\frac{\partial\;}{\partial z_{m}}
\end{align}
so that, along the leaf $L$, we have
\begin{align}\label{formulatwistleafhigher}  
{\mathcal L}^{j}_{m}= \sum_{l=n-p+1}^{n} \sum_{r}\left[(-\Lambda_{j r})\right]_{\bar m}{\cal H}^{r \bar l}\frac{\partial\;}{\partial z_{m}}
\end{align}
 \end{Prop}
\begin{pf} 
Computing Lie brackets using the frame ${\mathcal B}$ defined by (\ref{foliationframehigher}), the proposition follows immediately.
\end{pf}

\bigskip

As for codimension $1$ foliations, one may consider the Ricci curvature $\phi$ of the metric defined by
$dd^{c}u$ on the normal bundle to the foliation. 
Given a leaf $L$ and with respect to leaf coordinates $z_{1},\dots,z_{n}$ on $U$, so that $z_{1},\dots,z_{n-p}$ are holomorphic coordinates along  $L\cap U=\{z_{n-p+1}=\dots=z_{n}=0\}$ of the leaf, one has that $\phi$ is given by
\begin{align}\label{pformularicci}
\phi= - \frac{i}{2}\sum_{j,k=1}^{n-p}[\log (\det H)]_{j,\bar k}dz_{j}\wedge dz_{\bar k}.
\end{align}

Again, it is very important to remark that while the metric on the normal bundle to the foliation depends on the choice of the function $u$ the Ricci curvature does not  because of (\ref{annihilator}).
Exactly as  for codimension $1$ foliations one concludes:

\begin{Prop}  Let  \(\cal F\) be a codimension \(p\) locally Monge-Amp\`ere complex foliation  on a complex $n$-dimensional manifold $M$. 
The Ricci curvature locally defined by  (\ref{pformularicci}) defines globally  a Ricci curvature form of the normal bundle to the foliation  \(\cal F\) along any leaf of the foliation.
\end{Prop}

The next step is to provide a suitable expression for the Ricci curvature of the normal bundle to the foliation and to relate it to the twist tensor:

\begin{Prop} \label{ricciformhigher} If $z_{1},\dots,z_{n}$ are leaf coordinates  on an open set $U$, so that  $z_{1},\dots,z_{n-p}$  are holomorphic coordinates along the intersection $L\cap U=\{z_{n-p+1}=\dots=z_{n}=0\}$ of the leaf one has 
$$\phi = - \frac{i}{2}\sum_{j,k=1}^{n-p}S^{j\bar k}
dz_{j}\wedge dz_{\bar k}\hskip7.1cm$$
$$ =
- \frac{i}{2}\sum_{j,k=1}^{n-p}  \left(  \sum_{l,m,r,s=n-p+1}^{n} {\cal H}^{l\bar m}{\cal H}^{r \bar s}       
u_{\bar m j \bar s} u_{l\bar k r}  \right)dz_{j}\wedge dz_{\bar k}
$$
In particular for all $j=1,\dots,n-1$
$$0\leq S^{j\bar j}= \Vert {\cal L}^{j}\Vert ^{2}$$
where the norm is taken with respect to the metrics induced by $H$ and $H^{-1}={\cal H}$ on the normal bundle to the foliation and its dual respectively. Thus  the foliation is holomorphic if and only if these terms vanish.
 \end{Prop}
\begin{pf} As in the codimension $1$ case, the proof reduces to a computation which uses the formulas for the fourth derivatives  of the function $u$ along a leaf in terms of lower derivatives. Namely, with respect to leaf coordinates for a leaf  $L$,
 for $j, k\in \{1,\dots,n-p\}$ and indices $A,B\in \{1,\dots, n,\overline{1},\dots,\overline{n}\}$,
we have the following  along $L$:
\begin{align}\label{fourthhigher} u_{j\bar k AB}=  \sum_{l,m=n-p+1}^{n}
{\cal H}^{l\bar m}u_{l\bar k A}u_{j\bar m B}+  \sum_{r,s=n-p+1}^{n}
{\cal H}^{r \bar s}u_{r\bar k B}u_{j\bar s A}.
\end{align}
Moreover it is necessary that 
$A,B\in \{n-p+1,\dots, n,\overline{n-p+1},\dots,\overline{n}\}$ in order to have $u_{j\bar k AB}\neq0$ .
The proof of (\ref{fourthhigher}) is a similar to the one of Lemma \ref{forthderivatives} using this  time
the fact that any $(p+1)\times (p+1)$ minor containing $H$ has vanishing determinant and differentiating the equality. The computation is more tedious but completely elementary and it is very similar to the one carried out in \cite{Bedford-Burns}. The proof of the Proposition is then a repetition of the arguments given in Section 4 of \cite{Bedford-Burns} making use of (\ref{fourthhigher}).
\end{pf}

Suppose that ${\cal F}$ is a foliation of complex codimension $p$ on the complex $n$-dimensional manifold M locally defined in a neighborhood
U of $q \in M$ by the annihilator  of the form $dd^{c}u$ where $u$ is a plurisubharmonic
function on U. 
Let $L$ be a leaf of the Monge-Amp\`ere foliation. We need to show that all components of the the 
twist tensor ${\cal L}$ vanishes at all points $q\in L$. Since, by hypothesis $L$ is parabolic,  there is a holomorphic covering map $F\colon \C^{n-p} \to L$. For any $q\in L$, up to reparametrization, one may assume that $q=F(0)$. Let $I_{j}\colon \C\to \C^{n-1}$ be the $j$-th canonical injection: $I_{j}(z)=(0,\dots,z,\dots,0)$ and $f_{j}=F\circ I_{j} \colon \C\to  L$.
Then,  $L^{j}=f_{j}(\C)$ is a parabolic curve in $M$   for all $j=1,\dots,n-1$ with $q=f^{j}(0)\in L^{j}$. For all $j$ let $\iota_{j}\colon L^{j}\to L$  be the injection and denote with $\phi^{j} = \iota_{j}^{*}\phi$ the pull-back of the Ricci form $\phi$ of the metric on the normal bundle to $L$ defined by the form $dd^{c}u$. Then the form $\psi^{j}=-\phi^{j}$ is non-negative on $L^{j}$ and hence it defines a (pseudo$-$)metric on it. 

The key fact to get the main result on codimension $p$ foliations is the following:
\begin{Prop} 
With the notation as above, if $\psi^{j}(q)>0$, then at $q$, for each $j$
\begin{align}\label{riccihigher}
Ric(\psi^{j})\leq-\frac{2}{p}\psi^{j} 
\end{align}
\end{Prop}

\begin{pf} The proof is similar to the codimension $1$ case and it follows closely the proof of Theorem $3.1$ of
\cite{Burns}. Let $L$ be a leaf of the foliation and $q\in L$ a point where $\psi^{j}(q)>0$. Let
 $F\colon \C^{n-p} \to L$ be a holomorphic covering map, then $dF(0)\neq0$ and therefore 
$F$ provides local coordinates $z_{1},\dots,z_{n-p}$ along the leaf $L$ in a neighborhood (on $L$) of  $q$. We extend these coordinates to holomorphic coordinates $z_{1},\dots,z_{n-p}, \dots,z_{n}$ for $M$ in a neighborhood of $q$ so that, locally, $L$ is given by $\{z_{n-p+1}=\dots,z_{n}=0\}$, i.e. leaf coordinates. With respect to these coordinates, the point $q$ is the origin. Furthermore, 
 the computations of Proposition \ref{ricciformhigher} show that 
 $$ \psi^{j}= \frac{i}{2} S^{j\bar j}dz_{j}\wedge dz_{\bar j}= \left(
 \sum_{l,m,r,s=n-p+1}^{n} {\cal H}^{l\bar m}{\cal H}^{r \bar s}       
u_{\bar m j \bar s} u_{l\bar j r} \right) dz_{j}\wedge dz_{\bar j}
$$

For simplicity, we denote 
$$S= S^{j\bar j}=  \sum_{l,m,r,s} {\cal H}^{l\bar m}{\cal H}^{r \bar s}       
u_{\bar m j \bar s} u_{l\bar j r}.$$
and, as done in the last equality, from now on we shall just indicate  the indices in the sums with the understanding that all of them will range between $n-p+1$ and $n$. 
Proving our claim is equivalent to show the following equality:
\begin{align}\label{inequality} \left[\log S\right]_{j\bar j}(0)\geq \frac{2}{p} S(0).
\end{align}
Taking the $j$-th derivative of $S$, we have:
\begin{align}\label{Sj}
\begin{array}{ll}
S_{j}&= \left[ \displaystyle{\sum_{l,m,r,s}}{\cal H}^{l\bar m}{\cal H}^{r \bar s}       
u_{\bar m j \bar s} u_{l\bar j r}   \right]_{j}  \\
&\\
&=-\displaystyle{\sum_{l,h,k,m,r,s}}{\cal H}^{l\bar h} u_{k\bar h j} {\cal H}^{k\bar m}     {\cal H}^{r \bar s}       
u_{\bar m j \bar s} u_{l\bar j r}   \\
& \\
& \hskip1cm
-\displaystyle{\sum_{l,h,k,m,r,s}}  {\cal H}^{l\bar m}   {\cal H}^{r \bar h} u_{k\bar h j} {\cal H}^{k \bar s} 
u_{\bar m j \bar s} u_{l\bar j r}  \\
&\\
& \hskip1cm +
 \displaystyle{\sum_{l,m,r,s}}{\cal H}^{l\bar m}{\cal H}^{r \bar s}       
u_{\bar m j \bar s j} u_{l\bar j r} +
\displaystyle{\sum_{l,m,r,s}} {\cal H}^{l\bar m}{\cal H}^{r \bar s}       
u_{\bar m j \bar s} u_{l\bar j r j}  \\
&
\\
&= \displaystyle{\sum_{l,m,r,s}}
 {\cal H}^{l\bar m}{\cal H}^{r \bar s}   u_{l\bar j r}    
u_{jj\bar m \bar s } 
\end{array}
\end{align}
where, to get the last equality, we used (\ref{fourthhigher}) with $k=j$ and $A=l, B=r$ to substitute 
for $u_{l\bar j r j} $ and made the appropriate cancellations. Differentiating again, we have
\begin{align}\label{Sjbarjhigher}
\begin{array}{ll}
S_{j\bar j}&= \left[\displaystyle{\sum_{l,m,r,s}}{\cal H}^{l\bar m}{\cal H}^{r \bar s}   u_{l\bar j r}    
u_{jj\bar m \bar s }   \right]_{\bar j} \\
&\\
&= -\displaystyle{\sum_{l,h,k,m,r,s}}{\cal H}^{l\bar h} u_{k\bar h \bar j} {\cal H}^{k\bar m}   {\cal H}^{r \bar s}   
u_{l\bar j r}    
u_{jj\bar m \bar s } \\
&\\
&\hskip1cm
-\displaystyle{\sum_{l,h,k,m,r,s}} {\cal H}^{l\bar m} {\cal H}^{r \bar h} u_{k\bar h \bar j} {\cal H}^{k \bar s}   u_{l\bar j r}    
u_{jj\bar m \bar s } 
\\
&\\
&\hskip1cm + \displaystyle{\sum_{l,m,r,s}}{\cal H}^{l\bar m}{\cal H}^{r \bar s}   u_{l\bar j r \bar j}    
u_{jj\bar m \bar s }
+\displaystyle{\sum_{l,m,r,s}}{\cal H}^{l\bar m}{\cal H}^{r \bar s}   u_{l\bar j r}    
u_{jj\bar m \bar s  \bar j}.
\end{array}
\end{align}
To proceed, we would like to eliminate the fifth derivative term $u_{jj\bar m \bar s  \bar j}= u_{j\bar  j\bar m \bar s  j}$.
In order to do so we differentiate formula (\ref{fourthhigher})  with respect to $j$ in the case $k=j$, $A= \bar m$ and $B= \bar s$,  we substitute the forth order terms involving $j$ and $\bar j$ derivatives using again (\ref{fourthhigher}), finally getting:
\begin{align} \label{fifthhigher}
\begin{array}{ll}
u_{jj\bar m \bar s  \bar j} 
& = -\displaystyle{\sum_{a,b,c,d}} {\cal H}^{a\bar b} u_{c\bar b  j} {\cal H}^{c\bar d}      
u_{j\bar d m}    
u_{\bar j a \bar s } 
+ \displaystyle{\sum_{a,b}} {\cal H}^{a\bar b} u_{j j\bar b \bar m}    
u_{\bar j a \bar s}    
\\&
\\ & \hskip0,3cm 
+\displaystyle{\sum_{a,b,c,d}} {\cal H}^{a\bar b} u_{j\bar b \bar m} {\cal H}^{c\bar d}      
u_{j\bar d a}    
u_{\bar j c \bar s } 
+ 
\displaystyle{\sum_{a,b,c,d}}{\cal H}^{a\bar b} u_{j\bar b \bar m} {\cal H}^{c\bar d}      
u_{j\bar d \bar s}    
u_{\bar j c  a } 
\\&
\\ & \hskip0,3cm 
-\displaystyle{\sum_{a,b,c,d}}{\cal H}^{a\bar b} u_{c\bar b  j} {\cal H}^{c\bar d}      
u_{j\bar d \bar s}    
u_{\bar j a \bar m } 
+  \displaystyle{\sum_{a,b}}{\cal H}^{a\bar b} u_{j j\bar b \bar s}    
u_{\bar j a \bar m}    
\\ 
&\\ 
 & \hskip0,3cm +
 \displaystyle{\sum_{a,b,c,d}} {\cal H}^{a\bar b} u_{j\bar b \bar s} {\cal H}^{c\bar d}      
u_{j\bar d a}    
u_{\bar j c \bar m } 
+ 
\displaystyle{\sum_{a,b,c,d}}{\cal H}^{a\bar b} u_{j\bar b \bar s} {\cal H}^{c\bar d}      
u_{j\bar d \bar m}    
u_{\bar j c  a }.
\end{array}
\end{align}
Using (\ref{fifthhigher}) in (\ref{Sjbarjhigher}) and making the suitable cancellations,
one finally gets:
\begin{align}\label{Sjbarjfinal}
\begin{array}{ll}
S_{j\bar j}&= \displaystyle{\sum_{l,m,r,s}}
{\cal H}^{l\bar m}{\cal H}^{r\bar s} u_{j j \bar m \bar s } u_{\bar j \bar j l r } 
\\
&\\
& \hskip1cm +2  \displaystyle{\sum_{l,m,r,s,a,b,c,d}}
{\cal H}^{l\bar m}
{\cal H}^{r\bar s}
{\cal H}^{a\bar b}
{\cal H}^{c\bar d} 
u_{j\bar m \bar s}    
u_{j\bar b \bar d} 
u_{\bar j l  a }   
u_{\bar j r  c }. 
\end{array}
\end{align}
To get (\ref{inequality}), it is not restrictive to assume that at $0$ one has the normalization
$u_{r \bar s}= \delta_{r \bar s}$. Thus, again with the understanding that all the sum below extend over all indicated indices  running from $n-p+1$ to $n$, we have:
\begin{align}\label{S(0)}
S(0)=\displaystyle{\sum_{m,s}} \left| u_{j \bar m \bar s}\right|^{2},
\end{align}

\begin{align}\label{Sj(0)}
 \left|S_{j}(0)\right|^{2}= \left|\displaystyle{\sum_{l,r}}  u_{\bar j l r}    u_{jj\bar l \bar r }  \right|^{2}
 \leq 
 \displaystyle{\sum_{l,r}}  \left|u_{\bar j l r}  \right|^{2}  
 \displaystyle{\sum_{l,r}} \left| u_{jj\bar l \bar r }  \right|^{2}
 = S(0)\displaystyle{\sum_{l,r}}   \left|u_{jj\bar l \bar r }  \right|^{2}
\end{align}
and
\begin{align}\label{Sjbarj(0)}
S_{j\bar j}(0)= 
 \displaystyle{\sum_{l,r}}  \left|u_{jj\bar l \bar r }  \right|^{2} 
 + 2 \displaystyle{\sum_{m,s,b,d}} u_{j\bar m \bar s} u_{\bar j m  b } u_{j\bar b \bar d}   u_{\bar j s  d }.
 \end{align}
 Thus, if $S(0)\neq0$, we have
 \begin{align}\label{almostinequality}
S(0)\left[\log S\right]_{j\bar j}(0)= S_{j\bar j}(0) - \frac{\left|S_{j}(0)\right|^{2}}{S(0)}=
2  \displaystyle{\sum_{m,s,b,d}}  u_{j\bar m \bar s} u_{\bar j m  b } u_{j\bar b \bar d}   u_{\bar j s  d }.
 \end{align}

If we denote by $A$ the symmetric matrix of order $p$ defined by
$A=\left ( u_{\bar j s  d } \right)$ for $s,d=n-p+1\dots,n$ and by $A^{*}$  the conjugate 
transpose of $A$, then, using Schwarz's inequality,
 \begin{align}\label{inequalityequiv}\begin{array}{ll}
 S(0)\left[\log S\right]_{j\bar j}(0)
 & \geq
2  \sum u_{j\bar m \bar s} u_{\bar j m  b } u_{j\bar b \bar d}   u_{\bar j s  d } 
= 2\,Trace (A^{*}AA^{*}A) \\ 
&\\
&  \geq \frac{2}{p}\,\left[ Trace (A^{*}A)\right]^{2} 
= \frac{2}{p} \left[  \displaystyle{\sum_{m,s}}  \left| u_{j \bar m \bar s}\right|^{2} \right]^{2} =
\frac{2}{p} \left[  S(0) \right]^{2}.
\end{array}
 \end{align} 
If $S(0)\neq0$, (\ref{inequalityequiv})
is equivalent to (\ref{inequality})  and the proof is complete.
\end{pf}

We can now wrap up our discussion:

\begin{Th} \label{finaltheorem}A locally Monge-Amp\`ere foliation of rank $n-p$, where $1\leq p<n$, with parabolic leaves on a $n$-dimensional complex manifold is holomorphic. 
 \end{Th}
 
\begin{pf} We use the same notation as above. In the same vein of the argument given in \cite{Burns} 
and similarly to the proof of Theorem \ref{codimensionone}, we shall use inequality  (\ref{inequality}) and an Ahlfors Lemma argument to show that the Bedford-Burns twist tensor vanishes along all leaves of the foliation. 
Let $L$ be a leaf of the Monge-Amp\'ere foliation. We need to show that  the 
twist tensor ${\cal L}$ vanishes at all points $q\in L$. Suppose $q\in L$ is any point such that twist tensor $\cal L$ does not vanish at $q$.

Since, by hypothesis $L$ is parabolic,  there is a holomorphic covering map $F\colon \C^{n-p} \to L$. For any $q\in L$, up to reparametrization, one may assume that $q=F(0)$. Let $I_{j}\colon \C\to \C^{n-1}$ be the $j$-th canonical injection: $I_{j}(z)=(0,\dots,z,\dots,0)$ and $f_{j}=F\circ I_{j} \colon \C\to  L$.
Then,  $L^{j}=f_{j}(\C)$ is a parabolic curve in $M$   for all $j=1,\dots,n-p$ with $q=f_{j}(0)\in L^{j}$. For all $j$ let $\iota_{j}\colon L^{j}\to L$ the injection and denote with $\phi^{j} = \iota_{j}^{*}\phi$ the pull-back of the Ricci form $\phi$ of the metric on the normal bundle to $L$ defined by the form $dd^{c}u$. Then the form $\psi^{j}=-\phi^{j}$ is non negative on $L^{j}$ and hence it defines a (pseudo$-$)metric on it. 
For $R>0$, let 
$\D(R)=\{z\in \C \mid \vert z  \vert < R\}$. Then if we consider  $f_{j}\colon \D(R) \to L^{j}$, i.e. the restriction of the map $f_{j}$ to the disk $\D(R)$, then $\omega= {f_{j}}^{*}(\frac{2}{p}\psi^{j})$ satisfies 
$Ric(\omega)\leq -\frac{2}{p}\omega$ because of (\ref{inequality}). By Ahlfors' Lemma, then it follows that $\omega$ is dominated by the hyperbolic metric of $\D(R)$ defined by
the form
$$\omega_{0}= \frac{i}{2} \frac{R^{2}}{ (R^{2}-\vert z\vert^{2})^{2}}$$  
i.e. $\omega\leq \omega_{0}$ which, computing at the origin implies
\begin{align}\label{estimate} S(f(0))\vert {f_{j}}'(0)\vert^{2}\leq \frac{1}{R^{2}}.\end{align}
As $f_{j}$ is a covering map, ${f_{j}}'(0)\neq0$ so that, since (\ref{estimate}) holds for all $R>0$,
it follows that there is a contraddiction if $S(f^{j}(0))=S^{j\bar j}(q)\neq0$. 
In all these consideration $j=1,\dots,n-p$ was arbitrary and 
$q$ was any point in the leaf $L$.  By Proposition  3.2 then it follows that the twist tensor $\cal L$ vanishes at $q$, any $q\in L$. Since $L$ is any leaf, the foliation is holomorphic.
\end{pf}

\bigskip
\noindent
{\bf Remark.} For foliations of codimension $p>1$ the local Monge-Amp\`ere assumption we made in Theorem \ref{finaltheorem} cannot be removed. In fact in \cite{Duchamp-Kalka} it is discussed the codimension $2$ foliation of ${\bf C}^3$ in parabolic curves given by the affine part of an example due E. Calabi \cite{Calabi} of a foliation of ${\bf CP}^3$ by ${\bf C P}^1$. This foliation is neither holomorphic nor locally Monge-Amp\`ere. Namely, in \cite{Duchamp-Kalka} the Bedford-Burns twist tensor is explicitly computed, and is shown to be not vanishing and, furthermore, it is proved  that the foliation does not satisfy a symmetry condition which necessarily holds for locally Monge-Amp\`ere foliations.
See  \cite{Duchamp-Kalka} for the details.

\bigskip
\bigskip
\font\smallsmc = cmcsc8
\font\smalltt = cmtt8
\font\smallit = cmti8
\hbox{\parindent=0pt\parskip=0pt
\vbox{\baselineskip 9.5 pt \hsize=3.1truein
\obeylines
{\smallsmc
Morris Kalka
Mathematics Department
Tulane University
6823 St. Charles Ave.
New Orleans, LA 70118
USA
}\medskip
{\smallit E-mail}\/: {\smalltt kalka@@math.tulane.edu 
}
}
\hskip 0.0truecm
\vbox{\baselineskip 9.5 pt \hsize=3.7truein
\obeylines
{\smallsmc
Giorgio Patrizio
Dip. Matematica e Informatica``U. Dini''
Universit\`a di Firenze
Viale Morgani 67/a
I-50134 Firenze
ITALY
}\medskip
{\smallit E-mail}\/: {\smalltt patrizio@@unifi.it}
}
}

\end{document}